\documentclass[a4paper,12pt]{article} 

\usepackage{color}

\usepackage{amsfonts, amsmath, amsthm, amssymb}
\usepackage[T1]{fontenc}
\usepackage[cp1250]{inputenc}
\usepackage{xcolor}
\usepackage{graphicx}
\usepackage{amssymb}
\usepackage{amsmath}
\usepackage{mathptmx}
\usepackage{helvet}
\usepackage{courier}
\usepackage{txfonts}
\usepackage{tikz} 
\usetikzlibrary{arrows}
\usepackage{type1cm}

\usepackage{verbatim}

\usepackage{graphicx}
\usepackage{epsfig,amscd,amssymb,amsxtra,amsmath,amsthm}
\usepackage{type1cm}
\usepackage[T1]{fontenc}
\usepackage{graphics}
\usepackage[mathscr]{eucal}
\usepackage[all]{xy}
\usepackage{amsmath,amscd}

\newtheorem{theorem}{Theorem}[section]

\newtheorem{definition}[theorem]{Definition}

\newtheorem{observation}[theorem]{Observation}

\newcommand{\Cl}  {\mathop{\rm Cl}\nolimits}

\begin{document}

\def\joinrel{\mkern-3mu}
\newcommand{\varproj}{\displaystyle \lim_{\multimapinv\joinrel-\joinrel-}}

\title{An embedding of the Cantor fan into the Lelek fan}
\author{Iztok Bani\v c,  Goran Erceg,  and Judy Kennedy}
\date{}

\maketitle

\begin{abstract}
The Lelek fan $L$ is usually constructed as a subcontinuum of  the Cantor fan in such a way that the set of the end-points of $L$ is dense in $L$.   It easily follows that the Lelek fan is embeddable into the Cantor fan.  {It is also a well-known fact that the Cantor fan is embeddable into the Lelek fan, but this is less obvious.  When proving this, one usually uses the well-known result by Dijkstra and van Mill that the Cantor set is embeddable into the complete Erd\"os space,  and the well-known fact by Kawamura,  Oversteegen,  and Tymchatyn that the set of end-points of the Lelek fan is homeomorphic to the complete Erd\"os space.  Then, the subcontinuum  of the Lelek fan that is induced by the embedded Cantor set into the set of end-points of the Lelek fan,  is a Cantor fan.  

In our paper,  we give an alternative straightforward construction of a Cantor fan into the Lelek fan. We do not use the fact that the Cantor set is embeddable into the complete Erd\"os space and that it is homeomorphic to the set of end-points of the Lelek fan.  Instead, we use } our recent techniques of  Mahavier products of closed relations {to  produce an } embedding of the Cantor fan into the Lelek fan.  { Since the Cantor fan is universal for the family of all smooth fans,  it follows that also the Lelek fan is universal for smooth fans. }
\end{abstract}
\-
\\
\noindent
{\it Keywords:} Closed relations; Mahavier products;  fans; Cantor fans\;  Lelek fans\
\noindent
{\it 2020 Mathematics Subject Classification:} \emph{ 37B02, 37B45, 54C60, 54F15, 54F17}

\section{Introduction}

\begin{figure}[h!]
	\centering
		\includegraphics[width=16em]{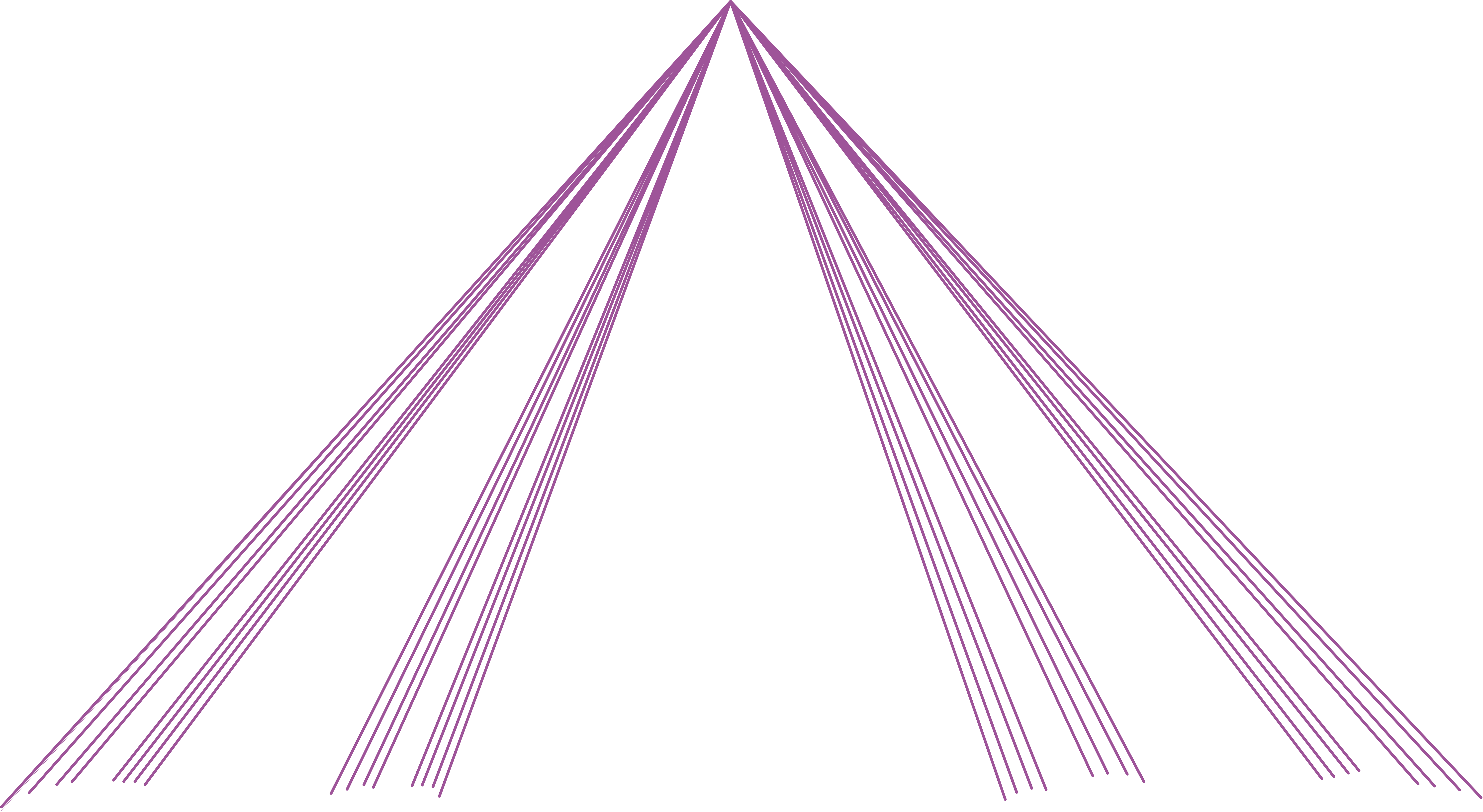}
	\caption{The Cantor fan}
	\label{fig000}
\end{figure}

 \emph{A continuum} is a non-empty compact connected metric space.  \emph{A subcontinuum} is a subspace of a continuum, which is itself a continuum.  Let $X$ be a continuum.  We say that $X$ is \emph{a Cantor fan}, if $X$ is homeomorphic to the continuum $\bigcup_{c\in C}A_c$, where $C\subseteq [0,1]$ is a Cantor set and for each $c\in C$, $A_c$ is the  {convex} segment in the plane from $(0,0)$ to $(c,-1)$; see Figure \ref{fig000}. 
Let $X$ be a Cantor fan and let $Y$ be a subcontinuum of $X$.  
A point $x\in Y$ is called an \emph{end-point of the continuum $Y$}, if for  every arc $A$ in $Y$ that contains $x$, $x$ is an end-point of $A$.  The set of all end-points of $Y$ will be denoted by $E(Y)$.   We say that the subcontinuum $Y$ of the Cantor fan $X$ is \emph{a Lelek fan}, if $\Cl(E(Y))=Y$.
\begin{figure}[h!]
	\centering
		\includegraphics[width=18em]{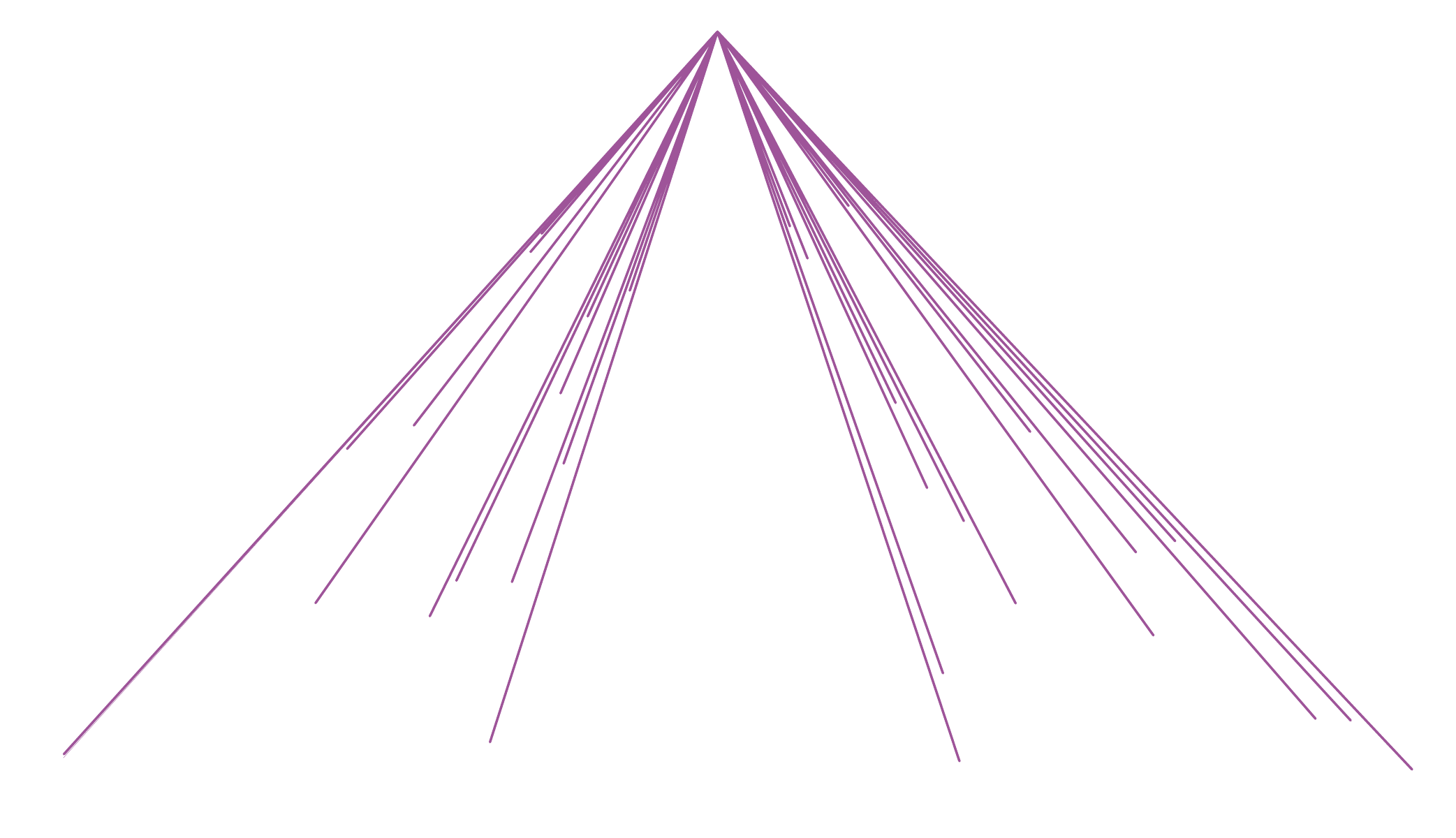}
	\caption{The Lelek fan}
	\label{figure2}
\end{figure}  
The first example of a Lelek fan was constructed by A.~Lelek in \cite{lelek}.  He proved  that the set of end-points of the Lelek fan is a dense one-dimensional set in the Lelek fan.  The Lelek fan  is also unique: any two non-degenerate subcontinua of the Cantor fan with a dense set of endpoints are homeomorphic.    This was proved independently by W.~D.~Bula and L.~Oversteegen  in \cite{oversteegen} and by W. ~Charatonik in  \cite{charatonik}.  After the Lelek construction, there were many other constructions of the Lelek fan. For example, in 2013,  D.~Bartosova and A.~Kwiatkowska constructed in \cite{dana} the Lelek fan as a quotient space of the projective Fraisse limit of a family that consists of finite rooted trees.  In \cite{banic1}, the Lelek fan is constructed  by I. Bani\v c,  G.  Erceg and J.  Kennedy as the inverse limit of inverse sequence of closed unit intervals with a single set-valued bonding function whose graph is an arc,  and in \cite{BEKMN}, the Lelek fan is constructed  by I. Bani\v c,  G.  Erceg, J.  Kennedy,  C.  Mouron and V.  Nall   as the inverse limit  of an inverse sequence of Cantor fans and a single transitive  continuous bonding function.

It easily follows from Lelek's construction  that the Lelek fan is embeddable into the Cantor fan.  However, it is  not that obvious that the Cantor fan is embeddable into the Lelek fan.  {One can easily construct an embedding of the Cantor fan into the Lelek fan by using
\begin{enumerate}
\item the well-known result  from \cite{dijkstra}  by J. ~J.~ Dijkstra and J. ~Mill that a space is almost zero-dimensional if and only if it is embeddable into the complete Erd\"os space, and 
\item the well-known result from \cite{kawamura}  by K. ~Kawamura, L. ~G.~Oversteegen, and E.~D.~Tymchatyn that the set of end-points of the Lelek fan is homeomorphic to the complete Erd\"os space.  
\end{enumerate} 
First, embed the Cantor set into the set of end-points of the Lelek fan and then,   the subcontinuum  of the Lelek fan that is induced by the embedded Cantor set,  is a Cantor fan (among other things, this was already noted by G.~Basso and R.~Camerlo in \cite{basso},  where another similar result is obtained). }

{ In this paper,  we give an alternative straightforward construction of a Cantor fan into the Lelek fan. In our approach, we do not use the well-known results from  \cite{dijkstra}  or \cite{kawamura}.   Instead, we use our recently developed techniques of  Mahavier products of closed relations from \cite{BEKMN},  \cite{banic1}, and \cite{banic2}. } We proceed as follows.  In Section \ref{s1}, the basic definitions and results that are needed later in the paper are presented. In Section \ref{s2}, our main result is proved.

\section{Definitions and Notation}\label{s1}
The following definitions, notation and well-known results will be needed in the paper.



\begin{definition}
Let $X$ be a {non-empty} compact metric space and let ${F}\subseteq X\times X$ be a relation on $X$. If ${F}$ is closed in $X\times X$, then we say that ${F}$ is  \emph{  a closed relation on $X$}.  
\end{definition}


\begin{definition}
Let $X$ be a {non-empty} compact metric space and let ${F}$ be a closed relation on $X$. Then we call
$$
X_F^+=\Big\{({x_0},x_1,x_2,x_3,\ldots )\in \prod_{k={0}}^{\infty}X \ | \ \textup{ for each {non-negative} integer } k, (x_{k},x_{k+1})\in {F}\Big\}
$$
\emph{ the  Mahavier product of ${F}$}.
\end{definition}

\begin{definition}
	For each $(r,\rho)\in (0,\infty)\times (0,\infty)$, we define the sets \emph{$L_r$}, \emph{$L_{\rho}$} and \emph{$L_{r,\rho}$}  as follows:
${L_r}=\{(x,y)\in [0,1]\times [0,1] \ | \ y=rx\}$,  ${L_{\rho}}=\{(x,y)\in [0,1]\times [0,1] \ | \ y=\rho x\}$, and $	{L_{r,\rho}}=L_r\cup L_{\rho}$.
	We also define the set \emph{$M_{r,\rho}$}  as follows:
	$$
	{M_{r,\rho}}{=[0,1]^+_{L_{r,\rho}}}. 
	$$
\end{definition}
\begin{definition}
	Let {$(r,\rho)\in (0,\infty)\times (0,\infty)$}. We say that \emph{$r$ and $\rho$ never connect} or \emph{$(r,\rho)\in \mathcal{NC}$}, if \begin{enumerate}
		\item $r<1$, $\rho>1$ and 
		\item for all integers $k$ and $\ell$,  
		$$
		r^k = \rho^{\ell} \Longleftrightarrow k=\ell=0.
		$$
	\end{enumerate} 
\end{definition}

%
In \emph{  \cite{banic1} }, the following theorem is the main result;  see 
\emph{ \cite[Theorem 14, page 21]{banic1}.}
\begin{theorem}\label{Lelek}
Let $(r,\rho)\in \mathcal{NC}$. Then $M_{r,\rho}$ is a Lelek fan with top $(0,0,0,\ldots)$.
\end{theorem}

In Theorem \ref{krajisca},  a characterization of end-points of $M_{r,\rho}$ is established; see \emph{ \cite[Theorem 3.5, page 8]{banic2}.  }

\emph{
\begin{definition}
	For each positive integer $k$,  we use $\pi_k:\prod_{i=1}^{\infty}[0,1]\rightarrow [0,1]$ to denote   the $k$-th standard projection from $\prod_{i=1}^{\infty}[0,1]$ to $[0,1]$.
\end{definition}}

\begin{theorem}\label{krajisca}
Let $(r,\rho)\in \mathcal{NC}$ and let $\mathbf x\in M_{r,\rho}$.  Then $\mathbf x\in E(M_{r,\rho})$ if and only if $\sup\{\pi_n(\mathbf x) \ | \ n \textup{ is a positive integer}\}=1$.
\end{theorem}

The following theorem is also  proved in \emph{  \cite[Theorem 9, page 18]{banic1}.}

\begin{theorem}\label{to}
Let $(r,\rho)\in \mathcal{NC}$. Then for each $x\in (0,1)$,  there is a sequence $a\in \{r,\rho\}^{\mathbb N}$ such that for each positive integer $n$,
$$
(a_1\cdot a_2\cdot a_3\cdot \ldots \cdot a_n)\cdot x \in [0,1]
$$
and 
$$
\sup\{(a_1\cdot a_2\cdot a_3\cdot \ldots \cdot a_n)\cdot x \ | \ n \textup{ is a positive integer}\}=1.
$$
\end{theorem}

\section{An embedding  of the Cantor fan  into the Lelek fan}\label{s2}
We show,  using our recent techniques from \cite{banic1} and \cite{banic2},  that the Cantor fan can be embedded into the Lelek fan.  

\begin{theorem}
The Cantor fan is embeddable into the Lelek fan.
\end{theorem}
\begin{proof}
Let $X=[0,1]$, let $(r,\rho)\in \mathcal{NC}$ and let 
$$
F=L_{r,\rho}\cup\{(t,t) \ | \  t\in [0,1]\} ~~~  \textup{ and }  ~~~ G=L_{r}\cup  \{(t,t) \ | \  t\in [0,1]\}.
$$
 It follows from \emph{\cite[Example 1, page 7]{banic1}} that $X_G^+$ is a Cantor fan.  Since $X_G^+\subseteq X_F^+$, it suffices to see that $X_F^+$ is a Lelek fan.  To do that, let 
 $$
 B_{\mathbf a}=\{(t,\mathbf a(1)\cdot t, \mathbf a(2)\mathbf a(1)\cdot t,\mathbf a(3)\mathbf a(2)\mathbf a(1)\cdot t,\ldots) \ | \ t\in [0,1]\}
 $$
 and
  $$
 A_{\mathbf a}=B_{\mathbf a}\cap X_F^+
 $$
 for each \emph{$\mathbf a = (\mathbf a(1), \mathbf a(2), \mathbf a(3), \ldots)) \in \{1,r,\rho\} ^ {\mathbb N}$.}
 Note that for each \emph{$\mathbf a\in \{1,r,\rho\}^ {\mathbb N}$}, $B_{\mathbf a}$ is a straight line segment in Hilbert cube $\prod_{k=1}^{\infty}[0,\rho^{k-1}]$ from $(0,0,0,\ldots)$ to $(1,\mathbf a(1)\cdot 1, \mathbf a(2)\mathbf a(1)\cdot 1,\mathbf a(3)\mathbf a(2)\mathbf a(1)\cdot 1,\ldots)$, and that for all \emph{$\mathbf a,\mathbf b\in \{1,r,\rho\}^ {\mathbb N}$}, 
 $$
 B_{\mathbf a}\cap B_{\mathbf b}=\{(0,0,0,\ldots)\} ~~~ \Longleftrightarrow ~~~ \mathbf a\neq \mathbf b.
 $$
Since \emph{$\{(1,\mathbf a(1)\cdot 1, \mathbf a(2)\mathbf a(1)\cdot 1,\mathbf a(3)\mathbf a(2)\mathbf a(1)\cdot 1,\ldots) \  | \  \mathbf a\in \{1,r,\rho\}^{\mathbb N}\}$} is a Cantor set, it follows that \emph{$\bigcup_{\mathbf a\in \{1,r,\rho\}^{\mathbb N}}B_{\mathbf{a}}$} is a Cantor fan.  Therefore,  $X_F^+$ is a  subcontinuum of the Cantor fan \emph{$\bigcup_{\mathbf a\in \{1,r,\rho\}^{\mathbb N}}B_{\mathbf{a}}$}.  Note that for each \emph{$\mathbf a\in \{1,r,\rho\}^{\mathbb N}$}, $A_{\mathbf a}$ is either degenerate or it is an arc from $(0,0,0,\ldots)$ to some other point, denote it by $\mathbf e_{\mathbf a}$.  Let 
\emph{$$
\mathcal U=\{\mathbf{a}\in  \{1,r,\rho\}^{\mathbb N} \ | \  A_{\mathbf a} \textup{ is an arc}\}.
$$}
Then 
$$
X_{F}^+=\bigcup_{\mathbf a\in \mathcal U}A_{\mathbf a}  ~~~ 
\textup{ and } ~~~ 
E(X_F^+)=\{\mathbf e_{\mathbf a} \ | \  \mathbf a\in \mathcal U\}.
$$
Next, we show that for each $\mathbf x\in X_F^+$,
$$
\mathbf x\in E(X_F^+) ~~~ \Longleftrightarrow ~~~ \sup\{\mathbf x(k) \ | \ k \textup{ is a positive integer}\}=1.
$$
Let $\mathbf x\in X_F^+$.  We treat the following possible cases. 
\begin{enumerate}
\item[Case 1.] For each positive integer $k$, there is a positive integer $\ell$ such that $\ell >k$ and $\mathbf x(k)\neq \mathbf x(\ell)$.  Then $\mathbf x\in E(X_F^+) ~~~ \Longleftrightarrow ~~~ \sup\{\mathbf x(k) \ | \ k \textup{ is a positive integer}\}=1$ follows from Theorem \ref{krajisca}.
\item[Case 2.] There is a positive integer $k$ such that for each positive integer $\ell\geq k$, $\mathbf x (\ell)=\mathbf x (k)$.  In this case, 
$$
\sup\{\mathbf x(k) \ | \ k \textup{ is a positive integer}\}=\max\{\mathbf x(k) \ | \ k \textup{ is a positive integer}\}
$$
and, therefore,  $\mathbf x\in E(X_F^+) ~~~ \Longleftrightarrow ~~~ \sup\{\mathbf x(k) \ | \ k \textup{ is a positive integer}\}=1$ easily follows. 
\end{enumerate}
To see that $X_F^+$ is a Lelek fan, let $\mathbf x\in X_F^+$ be any point and let $\varepsilon >0$. We prove that there is a point $\mathbf e\in E(X_F^+)$ such that $\mathbf e\in B(\mathbf x,\varepsilon)$ by considering the following possible cases.  
\begin{enumerate}
\item[Case 1.] For each positive integer $k$, there is a positive integer $\ell$ such that $\ell >k$ and $\mathbf x(k)\neq \mathbf x(\ell)$.   In this case,  the proof that there is a point $\mathbf e\in E(X_F^+)$ such that $\mathbf e\in B(\mathbf x,\varepsilon)$ is essentially the same as the proof  that there is a point $\mathbf e\in E(M_{r,\rho})$ such that $\mathbf e\in B(\mathbf x,\varepsilon)$ in the proof of Theorem \ref{Lelek}. Therefore, we leave the details to the reader.
\item[Case 2.] There is a positive integer $k$ such that for each positive integer $\ell\geq k$, $\mathbf x (\ell)=\mathbf x (k)$.  Without loss of generality, we assume that $\mathbf x\neq (0,0,0,\ldots)$.  Let $k_0$ be a positive integer such that $\sum_{k=k_0}^{\infty}\frac{1}{2^k}<\varepsilon$ and such that for each positive integer \emph{$k \geq k_0$, $\mathbf x (k)=\mathbf x (k_0)$}.  It follows from Theorem \ref{to} that there is a sequence $(a_1,a_2,a_3,\ldots)\in \{r,\rho\}^{\mathbf N}$ such that 
$$
\sup\{(a_1\cdot a_2\cdot a_3\cdot \ldots \cdot a_n)\cdot \mathbf x(k_0) \ | \ n \textup{ is a positive integer}\}=1.
$$
Choose and fix such a sequence $(a_1,a_2,a_3,\ldots)$.  Let
$$
\mathbf e=(\mathbf x(1), \mathbf x(2), \mathbf x(3), \ldots ,\mathbf x(k_0), a_1\cdot \mathbf x(k_0),a_2a_1\cdot \mathbf x(k_0),a_3a_2a_1\cdot \mathbf x(k_0),\ldots).
$$
Then $\mathbf e\in E(X_F^+)$ and 
$$
D(\mathbf e,\mathbf x)\leq \sum_{k=k_0}^{\infty}\frac{1}{2^k}<\varepsilon,
$$
where $D$ is the metric on $X_F^+$.
\end{enumerate}
This proves that $X_F^+$ is a Lelek fan. 
\end{proof}
\begin{observation} 
It is a well-known fact that the Cantor fan is universal for smooth fans (for details see \cite[Theorem 9, p. 27]{Jcharatonik},  \cite[ Corollary 4]{koch},  and \cite{eberhart}). Since the Lelek fan contains a Cantor fan, it follows  also that the Lelek fan is a universal continuum for smooth fans.  
\end{observation}

\section{Acknowledgement}
This work is supported in part by the Slovenian Research Agency (research projects J1-4632, BI-US/22-24-086, BI-HR/23-24-011 and BI-US/22-24-094, and research program P1-0285).
	

\noindent I. Bani\v c\\
              (1) Faculty of Natural Sciences and Mathematics, University of Maribor, Koro\v{s}ka 160, SI-2000 Maribor,
   Slovenia; \\(2) Institute of Mathematics, Physics and Mechanics, Jadranska 19, SI-1000 Ljubljana, 
   Slovenia; \\(3) Andrej Maru\v si\v c Institute, University of Primorska, Muzejski trg 2, SI-6000 Koper,
   Slovenia\\
             {iztok.banic@um.si}           
     
				\-
				
				\noindent G.  Erceg\\
             Faculty of Science, University of Split, Rudera Bo\v skovi\' ca 33, Split,  Croatia\\
{{goran.erceg@pmfst.hr}       }    

%
%
%
%
%
%
                 
                 	\-

  \noindent J.  Kennedy\\
             Lamar University, 200 Lucas Building, P.O. Box 10047, Beaumont, Texas 77710 USA\\
{{kennedy9905@gmail.com}       }


\end{document}